\documentclass[10pt,a4paper,oneside,leqno]{article}
\usepackage{amsfonts,amssymb}
\usepackage{amscd}
\newtheorem{defn}{Definition}

\newtheorem{observen}{Observation}

\vspace{0.5cm}
 4
\font\ebf=cmbx8
\font\erm=cmr8

\usepackage{graphicx, graphics}

\newcommand{\os}{\oplus\!\!\to}

\begin{document}
\begin{center}
	\noindent { \textsc{ Some Cobweb Posets Digraphs' Elementary Properties and Questions}}  \\ 
	\vspace{0.3cm}
	\vspace{0.3cm}
	\noindent Andrzej Krzysztof Kwa\'sniewski \\
	\vspace{0.2cm}
	\noindent {\erm Member of the Institute of Combinatorics and its Applications  }\\
{\erm High School of Mathematics and Applied Informatics} \\
	{\erm  Kamienna 17, PL-15-021 Bia\l ystok, Poland }\noindent\\
	\noindent {\erm e-mail: kwandr@gmail.com}\\
	\vspace{0.4cm}
\end{center}

\noindent {\ebf Abstract:}
\vspace{0.1cm}
\noindent {\small A digraph that represents reasonably  a scheduling problem should have no cycles i.e. it should be DAG i.e. a directed acyclic graph.  Here down we shall deal with special kind of graded DAGs named KoDAGs. For their definition and first primary properties see [1], where 
natural join of di-bigraphs (directed bi-parted graphs) and their corresponding  adjacency matrices is defined and then applied to investigate cobweb posets and their $Hasse$ digraphs called $KoDAGs$. 
\noindent In this report we extend the notion of cobweb poset while delivering  some elementary consequences of the description and observations established in [1].}

\vspace{0.3cm}


\vspace{0.1cm}

\noindent AMS Classification Numbers: 06A06 ,05B20, 05C7  
\vspace{0.1cm}

\noindent  affiliated to The Internet Gian-Carlo Polish Seminar:

\noindent \emph{http://ii.uwb.edu.pl/akk/sem/sem\_rota.htm}

\vspace{0.3cm}

\section{Introduction to the subject}

\vspace{0.2cm}
\noindent
It is now a Wiki important knowledge that an incidence structure is a triple $C = (P, L, I)$ where $P$ is a set of \emph{points}, $L$ is a set of \emph{lines} and  $I \subseteq P \times L$   is the incidence relation.  ( $I = P \times L$  for  KoDAGs  )
( compare:  $V = P \cup L, P \cap L = \emptyset$ ;  $P$ = black vertices = points, $L$= white vertices=lines).

\vspace{0.2cm}
\noindent
The elements of $I$ are called flags. If  $(p, l) \in I$  we say that point $p$ "lies on" line $l$.
The relation $I$  is equivalently  defined by its  bipartite digraph  $G( I )$. The relation $I$  and its bipartite digraph $G( I )$ are equivalently  defined by theirs  biadjacency matrix.  The example of thus efficiently coded finite geometries include such popular examples as Fano plane - a coding potrait of the distinguished composition algebra of  John T. Graves octonions (1843), a friend of William Hamilton,who called them octaves [2].   

\vspace{0.2cm}
\noindent
The \emph{incidence matrix} of an incidence structure $C$ is a \textbf{biadjacency matrix} of the Levi graph of the $C$ structure.

\vspace{0.2cm}
\noindent
The biadjacency matrix of a finite bipartite graph $G$ with n black vertices and $m$ white vertices is an $n \times m$ matrix where the entry $a_{ij}$ is the number of edges joining black vertex $i$ and white vertex $j$. In the special case of a finite, undirected, simple bipartite graph, the biadjacency matrix is a Boolean $(0,1)$-matrix.

\vspace{0.2cm}
\noindent
The adjacency matrix $A$ of a bipartite graph with the  reduced adjacency or -under synonymous substitution-  the biadjacency Boolean matrix $B$ is given by 

$$
	A = \left( \begin{array}{cc} 0 & B \\ B^T & 0 \end{array} \right).
$$

\vspace{0.2cm}
\noindent
The adjacency matrix $A$ of a bipartite \textbf{di}graph $\stackrel{\rightarrow}{K_{k,l}}$ (see: [1]) coded via its reduced adjacency or biadjacency Boolean matrix $B$ is according to [1] defined by

$$
	A = \left( 
	\begin{array}{cc} 0_{k,k} & B(k\times l) \\ 0_{l,k} & 0_{l,l} \end{array}
	\right),
	\quad where\ k=|P|,\ l=|L|
$$

\vspace{0.2cm}
\noindent \textbf{Example 1}

\begin{figure}[ht]
\begin{center}
	\includegraphics[width=70mm]{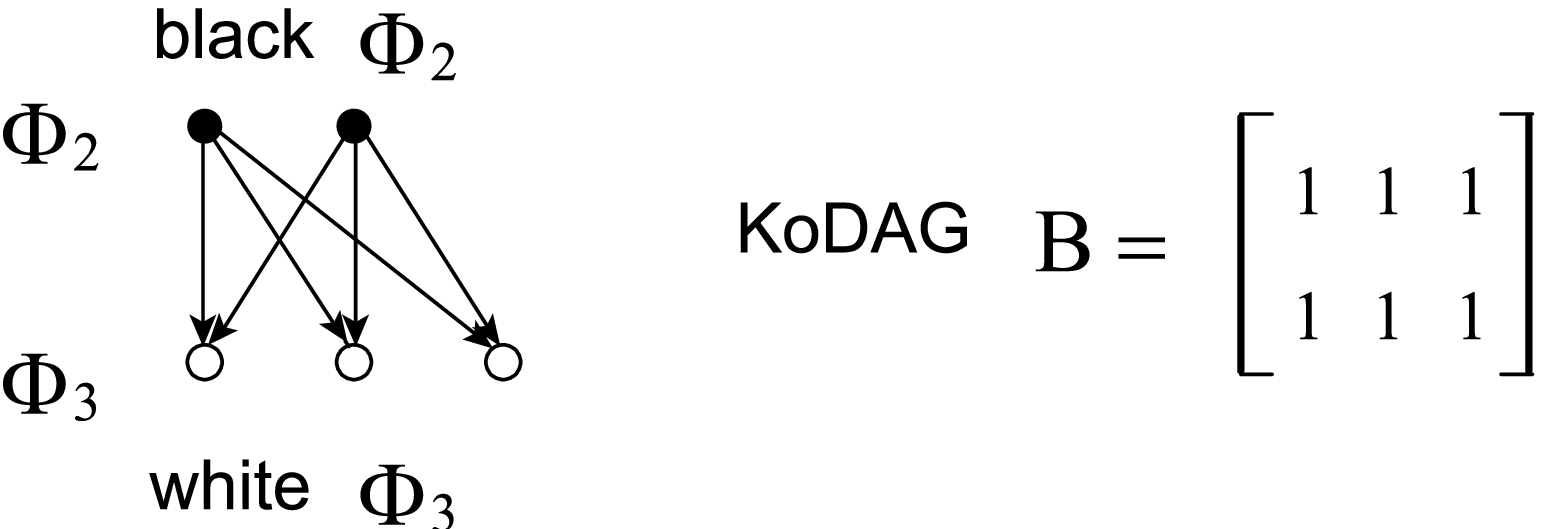}
	\caption{ \label{fig:1}}
\end{center}
\end{figure}

\begin{figure}[ht]
\begin{center}
	\includegraphics[width=70mm]{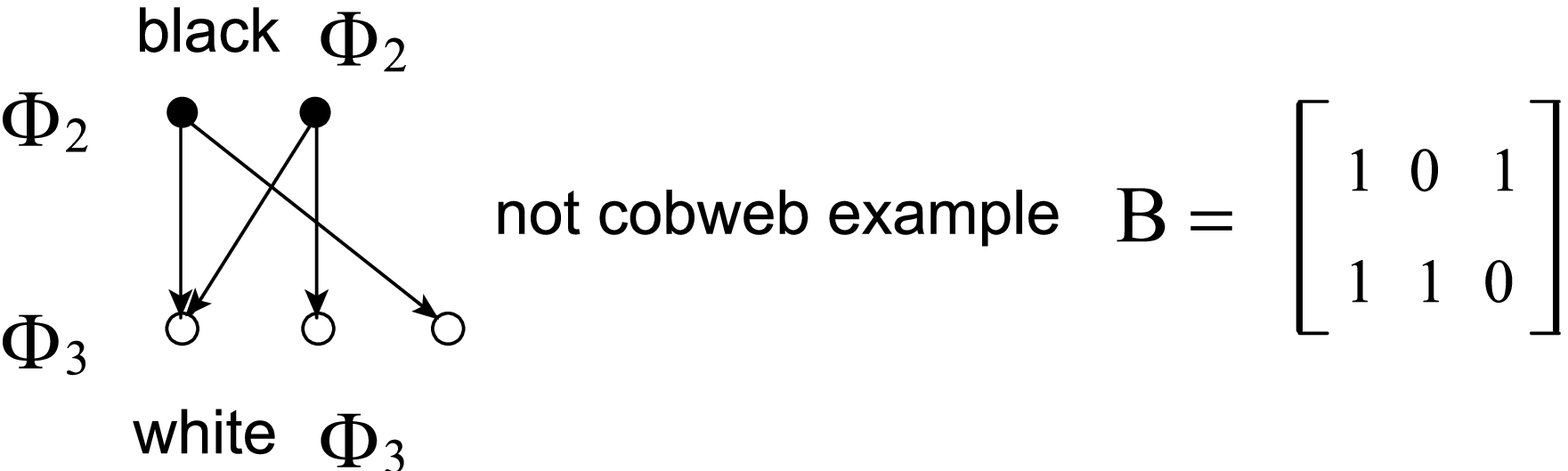}
	\caption{ \label{fig:1}}
\end{center}
\end{figure}

\vspace{0.2cm}
\noindent  Fig.1 displays (upside down way with respect to drawings in [1])  the bipartite digraph $\stackrel{\rightarrow}{K_{2,3}}$. It is obviously Ferrers dim 1 digraph [1].
Fig.2 displays   the bipartite sub-digraph of the $K$-digraph  $\stackrel{\rightarrow}{K_{2,3}}$. It is not Ferrers dim 1 digraph. What is its Ferrers dimension? Adjoin  minimal number of arcs in the Fig.2 in order to get  Ferrers dim 1 digraph,  bi-partite, of course. \\
The adjacency matrices coding digraphs from the  example above are shown below.

$$
	A_{KoDAG} = \left[ \begin{array}{cc}
	O_{2\times 2} & I(2\times 3) \\
	O_{3\times 2} & O_{3\times 3}
\end{array}\right] \ \ i.e. \ \
	A_{KoDAG} = \left( \begin{array}{cc}
	0_{2,2} & {1 1 1 \atop 1 1 1} \\
	0_{3,3} & 0_{3,3}
\end{array} \right),
$$
$$
	A_{not-cobweb} = \left( \begin{array}{cc}
	0_{2,2} & {1 0 1 \atop 1 1 0} \\
	0_{3,3} & 0_{3,3}
\end{array}\right)
$$

\noindent where  $O_{s \times s}$ stays for $(k \times m)$ zero matrix while  $I (s\times k)$  stays for  $(s\times k)$  matrix  of  ones  i.e.  $[ I (s\times k) ]_{ij} = 1$;  $1 \leq i \leq s$,  $1 \leq j \leq k$.  

\vspace{0.2cm}
\noindent Here above in the \textbf{Example 1} we are led  implicitly to the  notion of an extended cobweb poset as compared to [1] and references therein. For associated poset - see [1].

\begin{defn} ( extended cobweb poset- naturally graded)\\
Let $D = (\Phi, \prec\!\!\cdot )$  be a  transitive irreducible  digraph. Let  $n \in N \cup \{\infty\}$.  Let $D$ be a natural join $D = \os_{k=0}^n B_k$   of   Ferrers dim 1 bi-partite subdigraphs $B_k$ of di-bicliques  $\stackrel{\rightarrow}{K_{k,k+1}} = (\Phi_k \cup \Phi_{k+1}, \Phi_k\times\Phi_{k+1} ) , n \in N \cup \{\infty\}$.    The poset $\Pi (D)$ associated to this graded digraph $D = (\Phi,\prec\!\!\cdot )$    is called the extended cobweb poset or just cobweb, as a colloquial  abbreviation .
\end{defn}

\vspace{0.1cm}
\noindent Sometimes when we are in need we shall distinguish by name the \textbf{complete} cobwebs ( i.e. cobwebs represented by KoDAGs) from the overall  family of cobwebs (the extended cobweb posets as introduced above).

\vspace{0.2cm}

\noindent \textbf{Colligate with  Levi graph of an incidence structure}. Each incidence structure $C$ corresponds to a bipartite graph called Levi graph or incidence graph with a given black and white vertex coloring where black vertices correspond to points and white vertices correspond to lines of $C$ and the edges correspond to flags.

\vspace{0.2cm}
\noindent \textbf{Question 1}

\noindent
Is the natural join operation technique as started in [1] applicable to sequences of   Levi graphs of an incidence structures,somehow ?

\vspace{0.2cm}
\noindent In the case of graded digraphs with the finite set of minimal elements we have  what follows ( Observation 7  in [1]).

\begin{observen}
Consider bipartite digraphs'  chain obtained from the  di-biqliqes' chain via  deleting or no  arcs making thus [if deleting arcs] some or all of the di-bicliques $\stackrel{\rightarrow}{K_{k,k+1}}$  not di-biqliques; denote  them as  $G_k$. Let $B_k = B(G_k)$ denotes their biadjacency Boolean matrices correspondingly.  Then for any such $F$-denominated chain [hence any chain ] of bipartite digraphs  $G_k$  the general formula is:

$$
 B\left( \os_{i=1}^n G_i \right) \equiv  B [\os_{i=1}^n  A(G_i)] =  \oplus_{i=1}^n  B[A(G_i) ]  \equiv  \mathrm{diag} (B_1 , B_2 , ..., B_n) =
$$
$$
	= \left[ \begin{array}{lllll}
	B_1 \\
	& B_2 \\
	& & B_3 \\
	& & ... \\
	& & & & B_n
	\end{array} \right]
$$

\noindent $n \in N \cup \{\infty\}$.
\end{observen}

\vspace{0.2cm}
\noindent \textbf{Comment 1} (not only notation matter)

\noindent Let us denote by  $\langle\Phi_k\to\Phi_{k+1}\rangle$  the di-bicliques  denominated by subsequent levels $\Phi_k, \Phi_{k+1}$ of the graded  $F$-poset $P(D) = (\Phi, \leq)$  i.e. levels $\Phi_k , \Phi_{k+1}$ of  its cover relation graded digraph  $D = (\Phi,\prec\!\!\cdot$)  i.e. Hasse diagram  (see notation in the authors and others papers quoted in [1]). Then one  may conditionally  approve the following identification if necessary natural join condition [1] is implicit within this identification. 

$$
	B\left(\os_{k=1}^n \langle\Phi_k\to\Phi_{k+1}\rangle \right) = B\left(\bigcup_{k=1}^n \langle\Phi_k\to\Phi_{k+1}\rangle \right)  $$
if the conditioned set sum of digrahps concerns an ordered digraphs's pair satisfying natural join condition [1] what makes such a conditioned  set sum of vertices and simultaneously the set sum of disjoint arcs $E_k$ , $E_{k+1}$ families non commutative.
\noindent  Note that this what  just has been  said is exactly the reason of  $B(G_1\os G_2) = B(G_1\cup G_2) = B(G_1)\oplus B(G_2)$.

\section{On number of finite cobwebs an related questions}
\subsection{Two schemes and a Question}
\noindent Before we deal with  questions "`how many"' let us jot first two schemes of two statements which may be simultaneously referred to
relations, their digraphs or corresponding adjacency matrices. Secondly comes an elementary question without giving an answer. 

\vspace{0.2cm}

\begin{center}
\textbf{(Ferrers dim 1)}$\os$ \textbf{(Ferrers dim 1)}   = \textbf{(Ferrers dim 1)}.
\end{center}

\vspace{0.1cm}
\noindent(Obvious: use  $2\times2$ permutation sub-matrix forbidding i.e.  $2\times2$ permutation sub-matrix disqualification criterion)

\vspace{0.2cm}

\begin{center}
\textbf{Ferrers}  $\os$\textbf{ Ferrers } = \textbf{Ferrers}.
\end{center}

\vspace{0.2cm}
\noindent  See Observation 3  in [1]  and note that resulting  biadjacency matrices neither  contain  any of  two $2\times2$ permutation matrices. 
Nota bene the Observation 3 from [1]  follows from the above obvious statements.

\vspace{0.2cm}
\noindent \textbf{Question 2}

\noindent For  biadjacency matrices  $B(G_1) = B_1$   and  $B(G_2) = B_2$ of bipartite digraphs $G_1$ and $G_2$  we have the matrix exponential rule

\vspace{0.1cm}

$$\exp[B_1\oplus B_2] =  \exp[B_1] \otimes \exp[B_2] \:,$$
where $\otimes$   stays for the Kronecker product.

\vspace{0.2cm}
\noindent Let $F$ be any natural number valued sequence. Let $\mathbf{A}_F$  denotes the Hasse matrix of  the $F$-denominated cobweb poset $\langle\Phi,\leq\rangle$ [1].   The $\zeta$   matrix is then the  geometrical series in  $\mathbf{A}_F$: $\zeta = (1 - \mathbf{A}_F)^{-1 \copyright}$.
(Recommended: consult the remark from page 12 in [1] on $\zeta = exp[\mathbf{A}]$ in the cases of the Boolean  poset $2^N$ and the Ferrand-Zeckendorf poset of finite subsets of $F=N$  without two consecutive elements [3]). The \textbf{Question 2} is : find the rule if any for 
  
$$\zeta_{B_1 \oplus B_2} = (1 - B_1\oplus B_2)^{-1 \copyright}= ?$$

\subsection{How many}

\noindent \textbf{Notation for this subsection.}

\noindent Consider natural number \textit{N} composition $N = f_1 + f_2 +...+ f_k$ where  $ 0< f_1,f_2,...,f_k \leq N$. The compositions' \textbf{type} 
$\stackrel{\rightarrow}{k} = \left\langle f_1,f_2,...,f_k \right\rangle$ labels compositions of the chosen natural number $N$, where $N=|V|$ labels on its own the partial graded orders  $P_N= \left\langle V, \leq \right\rangle$ with $N$ points (vertices) and the partition 
$V =  \bigcup_{r=1}^k V_r$, $V_r \cap V_s = \oslash $ for $r \neq s$,  $f_r= |V_r|$ and $r,s =1,...,k$ , $k=1,...,N$. The partial order $\leq $ is the subset according to $ \leq \subseteq V_1 \times V_2 \times ...\times V_k$.  The symbol $\Big\{{N \atop k}\Big\}$ denotes the array of Stirling numbers of the second  kind.

\vspace{0.2cm}
\noindent \textbf{Obvious from obvious and Questions}

\noindent Number of all $k$-tuples for any $k$-block ordered partition  $<V_1,V_2,…,V_k>$    equals to 

\begin{center}
              $|V_1|\cdot|V_2|\cdot ...\cdot |V_k|=  \prod_{r=1}^k V_r.$      
\end{center}

\vspace{0.2cm}
\noindent The number of all complete cobweb posets  $P_N= \left\langle V, \leq \right\rangle$  with  $|V|=N$  elements is equal
to  $T_N$ = the number of ordered partitions of V. - Why?   Note:    The number of ordered partitions of  $\left\langle f_1,f_2,...,f_k \right\rangle$   = $\stackrel{\rightarrow}{k}$ type is equal  to $ {n \choose {f_1,f_2,...,f_k}} = \frac{n!}{f_1!f_2!...f_k!} $.  Thereby:
the number of all complete cobweb posets  $P_N= \left\langle V, \leq \right\rangle$ of  $\left\langle f_1,f_2,...,f_k \right\rangle$= $\stackrel{\rightarrow}{k}$ type is equal  to   $ {n \choose {f_1,f_2,...,f_k}}$ where $f_r= |V_r|$ and $r =1,...,k$  for all particular $\stackrel{\rightarrow}{k}$ type $k$-block ordered partitions  $\bigcup_{r=1}^k V_r = V$.     \textbf{Altogether}:

\vspace{0.2cm}

\noindent \textbf{2.2.1.} The number $Cob^c(N,\stackrel{\rightarrow}{k})$  of all complete of the type  $\left\langle f_1,f_2,...,f_k \right\rangle \equiv \stackrel{\rightarrow}{k}$  cobweb posets $P_N$ is given by:

\begin{center}
$Cob^c(N,\stackrel{\rightarrow}{k})   =  {N\choose {f_1,f_2,...,f_k}}$   ,
  $N = f_1 + f_2 +...+ f_k$ ,  $ 0< f_1,f_2,...,f_k \leq N$,
\end{center}

\vspace{0.3cm}

\noindent \textbf{2.2.2.} The number $Cob^c(N,k)$  of all complete \textbf{k}-level cobweb posets $P_N$ reads

$$Cob^c(N,k) = \sum_{{f_1+f_2+...+f_k = N} \atop {0<f_1,...,f_k\leq N}}{N \choose {f_1,f_2,...,f_k}}= k!\Bigg\{ {N \atop k} \Bigg\} =
  \sum_{r=0}^k (-1)^{N-k} r^N {N \choose r}.$$

\begin{center}
                    $Cob^c(N,k)$ =  number of surjections   $f: V \mapsto  [k].$
\end{center}

\vspace{0.3cm}

\noindent \textbf{2.2.3.} The number  $Cob^c(N)$  of all complete cobweb posets $P_N:$ is then the sum:

 $$Cob^c(N)  =  \sum_{k=1}^N k!\Bigg\{ {N \atop k} \Bigg\}.$$ 
  
\begin{center}
$Cob^c(N) = T_N$ = the number of \textbf{ordered} partitions of $V$.
\end{center}

\vspace{0.3cm}

\noindent \textbf{2.2.4.} The number $K(N,\stackrel{\rightarrow}{k})$   of all \textbf{\textsl{k}}-ary relations  of the given $\stackrel{\rightarrow}{k}$ type is:

$$K(N,\stackrel{\rightarrow}{k}) = 2^{\prod_{r=1}^k V_r} - 1,$$
where 
\begin{center}
$N = f_1 + f_2 +...+ f_k$ ,$ 0< f_1,f_2,...,f_k \leq N$, $f_r= |V_r|$ and $r=1,...,k$ while  $k=1,...,N.$
\end{center}

\vspace{0.3cm}

\noindent \textbf{2.2.5.} For the number $K(N)$  of all type $k$-ary relations , $k = 1,…, N$ we then have

$$ K(N) = \sum_{{f_1+f_2+...+f_k = N} \atop {0<f_1,...,f_k\leq N}}[2^{f_1\cdot f_2\cdot...\cdot f_k} - 1].$$

\vspace{0.3cm}

\noindent \textbf{2.2.6.} \textbf{Question 3.}  The number of all  \textbf{\textit{k}} -level graded posets $P_N= \left\langle V, \leq \right\rangle$  
with $|V|=N$ elements where the partial order $\leq $ is the subset of Cartesian product: $ \leq \subseteq V_1 \times V_2 \times ...\times V_k$
and  where $f_r= |V_r|$ , $r =1,...,k$ and  $k=1,...,N$  while  $N = f_1 + f_2 +...+ f_k$ ,  $ 0< f_1,f_2,...,f_k \leq N$ ... \textbf{equals ?}

\vspace{0.3cm}

\noindent \textbf{2.2.7.} \textbf{Question 4.} The number of all graded posets $P_N= \left\langle V, \leq \right\rangle$ with  $|V|=N$  elements for any  type  $\stackrel{\rightarrow}{k}$, $k = 1,…, N$ ... \textbf{equals ?}

\vspace{0.4cm}

\noindent \textbf{Acknowledgments}.  Let me  thank the Student of Gda\'nsk University Maciej Dziemia\'nczuk for his alert assistance.  I am indebted  with gratitude to my wife for being so patient with me still browsing cobwebs in spite of rolling around Christmas Eve.

\end{document}